\documentclass[12pt, a4paper]{article}
\usepackage[latin1]{inputenc}
\usepackage[english]{babel}
\usepackage{indentfirst}
\usepackage{amstext}
\usepackage{amsfonts}
\usepackage{textcomp}
\usepackage{amssymb}
\usepackage[dvips]{graphicx}
\usepackage{setspace}
\usepackage[all]{xy}

\newcommand {\demo}{\hskip -0.6cm{\bf Proof.  }}
\newcommand {\fim}{\hfill{$\square$}\vskip 1pc}

\newcommand {\R}{\mathbb{R}}
\newcommand {\N}{\mathbb{N}}
\newcommand {\Z}{\mathbb{Z}}

\newcommand{\supp}{\text{supp}}

\newtheorem{teorema}{Theorem}[section]

\newtheorem{corolario}[teorema]{Corollary}
\newtheorem{definicao}[teorema]{Definition}
\newtheorem{proposicao}[teorema]{Proposition}
\newtheorem{exemplo}[teorema]{Example}

\begin{document}

\onehalfspace

\title{Graph C*-algebras, branching systems and the Perron-Frobenius operator}
\author{Daniel Gonçalves and Danilo Royer}
\date{}
\maketitle

AMS 2000 MSC: 47L99, 37A55


\begin{abstract}
In this paper we show how to produce a large number of representations of a graph C*-algebra in the space of the bounded linear operators in $L^2(X,\mu)$. These representations are very concrete and, in the case of graphs that satisfy condition $(L)$, we use our techniques to realize the associated graph C*-algebra as a subalgebra of the bounded operators in $L^2(\R)$. We also show how to describe the Perron-Frobenius operator of ergodic theory in terms of the representations we associate to a graph.
\end{abstract}

\section{Introduction}

The theory of graph algebras has been explored extensively in recent years, both in pure algebra and in operator theory, see e.g., \cite{GR3, KST09, flr, Po09}. Similarly branching systems arise in neighboring disciplines such as random walk, symbolic dynamics and scientific computing,
see e.g., \cite{BGT06,  Dev07, HG09, LM09, SaSt96}. 

Our aim in this paper is to explore the richness of the theories of graph C*-algebras and branching systems and relate then. The pillar of our work is the construction of representations of graph C*-algebras from branching systems. As a consequence of the study of these representations we obtain a concrete description of the Perron-Frobenius operator of ergodic theory, as well as faithfull representations of many graph C*-algebras (including the algebra of compact operators in a separable Hilbert space). 

We expect that our theory will have many more applications, much in the same manner as the theory of representations of the Cuntz algebras is crucial to the understanding of representations of the fermion algebra, see \cite{AK02}, the classification of theories of quantum string fields, see \cite{AK03}, and is also applied to dynamical systems, see \cite{GR, K20031, K20032}, fractals see \cite{K20033} and the theory of wavelets, see \cite{bratteli}. We expect that many of the results in the literature above can be generalized to the graph C*-algebra setting. For example, we have already explored some of these possibilities in \cite{GR2}.

Given a graph $E$, in \cite{flr}, the associated graph C*-algebra is defined as a universal C*-algebra generated by projections and partial isometries satisfying given relations. Even though this definition is completely clear, the use of a universal object brings a level of abstraction, which sometimes may elude the non-expert. As a consequence of our study we are able to give a concrete characterization of graph C*-algebras (for graphs that satisfy condition $(L)$) as subalgebras of the bounded operators in $L^2(\R)$, $B(L^2(\R))$. This includes many know algebras, as for example the algebra of compact operators, and hence, if the reader so desire, it could define the compact operators as a subalgebra generated by multiplication and composition operators in $B(L^2(\R))$ (see example \ref{compact}). For graphs in general, we show how to obtain representations of the associated graph C*-algebra in $B(L^2(\R))$, but, without the presence of condition $(L)$, we can not guarantee that these representations are faithful. Still, for any countable graph, we show how the representations mentioned above can be used to describe the Perron-Frobenius operator in $L^1(X,\mu)$ (this is the analogue of what was done in \cite{GR} and \cite{kawamura} for the algebra $O_A$), and so we establish a link between the operator theory of graph C*-algebras and the ergodic theory of nonsingular maps. 

The paper is divided in five sections. After this brief introduction, in section 2, we introduce E-branching systems associated to a graph and show how they induce representations of the graph C*-algebra. Next, in section 3, we prove the existence of E-branching systems in $\R$ for any graph with countable edges and vertices and show that, for graphs that satisfy condition $(L)$, the representations arising from these E-branching systems in $\R$ are faithful. In section 4, we show how the representations mentioned above relate with the Perron-Frobenius operator and we finish the paper in section 5, where we present two examples. Before we proceed, we recall the definition of graph C*-algebras below.

Let $E=(E^0, E^1, r, s)$ be a directed graph, that is, $E^0$ is a set of vertices, $E^1$ is a set of edges and $r,s:E^1\rightarrow E^0$ are the range and source maps. Following \cite{flr}, the C*-algebra of the graph $E$ is the universal C*-algebra, $C^*(E)$, generated by projections $\{P_v\}_{v\in E^0}$ and partial isometries $\{S_e\}_{e\in E^1}$ with orthogonal ranges satisfying:
\begin{itemize}
\item the projections $p_v$ are mutually orthogonal,
\item $S_e^*S_e=P_{r(e)}$ for each $e\in E^1$,
\item $S_eS_e^*\leq P_{s(e)}$ for each $e\in E^1$,
\item $P_v=\sum\limits_{e:s(e)=v}S_eS_e^*$ for every vertex $v$ with $0<\#\{e:s(e)=v\}<\infty$.
\end{itemize}

\section{E-branching system}

In this section we will define the E-branching system associated to a directed graph $E$ and we will show that each E-branching system induces a representation of the graph algebra $C^*(E)$. Before we proceed we would like to mention that even though the definition of an E-branching system seems rather technical, it is nothing more than the translation of the conditions in the definition of graph C*-algebras to the measurable setting, as we shall see below.

Throughout the paper we will use some notation about operations over measurable sets and maps. For measurable subsets $A,B$ in a given measure space $(X,\mu)$, the notation $B\stackrel{\mu-a.e.}{\subseteq}A$ means that $\mu(B\setminus A)=0$, and the notation $A\stackrel{\mu-a.e.}{=}B$ means that $\mu(A\setminus B)=0$ and $\mu(B\setminus A)=0$. For two maps, $f,g:A\rightarrow X$, the notation $f\stackrel{\mu-a.e.}{=}g$ means that $\mu(x\in A:f(x)\neq g(x))$=0.

\begin{definicao}\label{brancsystem}
Let $(X,\mu)$ be a measure space and let $\{R_e\}_{e\in E^1}$, $\{D_v\}_{v\in E^0}$ be families of measurable subsets of $X$ such that:
\begin{enumerate}
\item $R_e\cap R_d\stackrel{\mu-a.e.}{=} \emptyset$ for each $d,e\in E^1$ with $d\neq e$,
\item $D_u\cap D_v\stackrel{\mu-a.e.}{=} \emptyset$ for each $u,v\in E^0$ with $u\neq v$,
\item $R_e\stackrel{\mu-a.e.}{\subseteq}D_{s(e)}$ for each $e\in E^1$,
\item $D_v\stackrel{\mu-a.e.}{=} \bigcup\limits_{e:s(e)=v}R_e$\,\,\,\,\, if\,\,\,\,\, $0<\#\{e\in E^1\,\,:\,\,s(e)=v\}< \infty$,
\item for each $e\in E^1$, there exists a map $f_e:D_{r(e)}\rightarrow R_e$ such that $f_e(D_{r(e)})\stackrel{\mu-a.e.}{=}R_e$ and the Radon-Nikodym derivative $\Phi_{f_e}$ of $\mu\circ f_e$, with respect to $\mu$ (in $D_{r(e)}$), exists and $\Phi_{f_e}> 0$ $\mu$ a.e.,
\item for each $f_e$ as above there exists a map $f_e^{-1}:R_e \rightarrow D_{r(e)}$ such that $f_e\circ f_e^{-1}\stackrel{\mu-a.e.}{=}Id_{R_e}$ and  $f_e^{-1}\circ f_e\stackrel{\mu a.e.}{=}Id_{D_{r(e)}}$, and for each such $f_e^{-1}$ there exists the Radon-Nikodym derivative $\Phi_{f_e^{-1}}$ of $\mu\circ f_e^{-1}$ with respect to $\mu$ (in $R_e$).
\end{enumerate}

A measurable space $(X,\mu)$, with families of measurable subsets $\{R_e\}_{e\in E^1}$ and $\{D_v\}_{v\in E^0}$, and maps $f_e$, $f_e^{-1}$, $\Phi_{f_e}$ and  $\Phi_{f_e^{-1}}$ as above is called an $E$-branching system.
 
\end{definicao}

In the fifth item from the definition above, the domain of the measures $\mu\circ f_e$ and $\mu$ are the measurable subsets of $D_{r(e)}$. So the Radon-Nikodym derivative $\Phi_{f_e}$ is a measurable map with domain $D_{r(e)}$. We will consider $\Phi_{f_e}$ also as a measurable map with domain $X$ (defining it as being zero out of $D_{r(e)}$). The same holds for the map $\Phi_{f_e^{-1}}$, that is, $\Phi_{f_e^{-1}}$ will be considered as a measurable map with domain $R_e$ and $X$.

The next step is to show that each $E$-branching system induces a representation of the $C^*$-algebra $C^*(E)$. So, let $(X,\mu)$ be an $E$-branching system, as in the definition above. For each $e\in E^1$, define the operator $\pi(e)\in B(L^2(X,\mu))$ (the bounded linear operators in $L^2(X,\mu)$) as follows: for each $\phi\in L^2(X,\mu)$, and $x\in R_e$, let
$$\pi(e)\phi_{|_{x}}=\Phi_{f_e^{-1}}^{\frac{1}{2}}(x)\phi(f_e^{-1}(x))$$
and if $x\notin R_e$, let $\pi(e)\phi_{|_{x}}=0$. 

In order to simplify notation, in what follows we will make a small abuse of the characteristic function symbol and denote the above operator as: 
$$\pi(e)\phi=\chi_{R_e}\cdot \Phi_{f_e^{-1}}^{\frac{1}{2}}\cdot \phi\circ f_e^{-1}.$$ 
It is easy to show that $\pi(e)\phi\in L^2(X,\mu)$, for each $\phi\in L^2(X,\mu)$. Also, $\pi(e)$ is linear and $||\pi(e)\phi||\leq ||\phi||$, and so $\pi(e)$ lies in fact in $B(L^2(X,\mu))$. The adjoint of $\pi(e)$ is the operator defined as:
$$\pi(e)^*\phi=\chi_{D_{r(e)}}\cdot\Phi_{f_e}^{\frac{1}{2}} \cdot \phi \circ f_e,$$ (where we are using the characteristic function symbol with the same meaning as in the notation of $\pi(e))$.

For each $v\in E^0$, define $\pi(v):L^2(X,\mu)\rightarrow L^2(X,\mu)$ by 

$$\pi(v)\phi=\chi_{D_v}\phi, \text{ for all }\phi\in L^2(X,\mu).$$
 (that is, $\pi(v)$ is the multiplication operator by $\chi_{D_v}$, the characteristic function of $D_v$). 

 \begin{teorema}\label{rep} Let $(X,\mu)$ be an $E$-branching system. Then there exists a *-homomorphism $\pi:C^*(E)\rightarrow B(L^2(X,\mu))$ such that 
$$\pi(S_e)\phi=\chi_{R_e}\cdot\Phi_{f_e^{-1}}^{\frac{1}{2}}\cdot \phi\circ f_e^{-1} \text{ and }\pi(P_v)\phi=\chi_{D_v}\phi,$$ for each $e\in E^1$ and $v\in E^0$.
  \end{teorema}
  
\demo 
For each $e\in E^1$ and $v\in E^0$ define $\pi(S_e)=\pi(e)$ and $\pi(P_v)=\pi(v)$, where $\pi(e), \pi(v)$ are as above. Note that $\pi(S_e)^*\pi(S_e)=M_{D_{r(e)}}$ and $\pi(S_e)\pi(S_e)^*=M_{R_e}$. To obtain the desired *-homomorphism $\pi:C^*(E)\rightarrow B(L^2(X,\mu))$ it is enough to verify that the families $\{\pi(S_e)\}_{e\in E^1}$ and  $\{\pi(P_v)\}_{v\in E^0}$ satisfy the relations which define $C^*(E)$. Obviously $\pi(S_e)$ are partial isometries and $\pi(P_v)$ are projections, for all $e\in E^1$ and $v\in E^0$. Note that the projections $\pi(P_v)$ are mutually orthogonal, because $D_u\stackrel{\mu- a.e.}{\cap} D_v=\emptyset$ for $u\neq v$. The equality $\pi(S_e)^*\pi(S_e)=\pi(P_{r(e)})$ is immediate, and the inequality $\pi(S_e)\pi(S_e)^*\leq \pi(P_{s(e)})$  follows from the third item of the E-branching system definition. To verify the last relation, let $v\in E^0$ be such that $0<\{e\in E^1\,:\, s(e)=v\}<\infty$. Then $D_v=\bigcup\limits_{e:s(e)=v}R_e$, and so $$M_{\chi_{D_v}}=M_{\chi_{\bigcup\limits_{e:s(e)=v}R_e}}.$$ Since $R_e\stackrel{\mu-a.e.}{\cap}R_d=\emptyset$ for $e\neq d$, then  

$$M_{\chi_{\bigcup\limits_{e:s(e)=v}R_e}}=\sum\limits_{e:s(e)=v}M_{\chi_{R_e}}.$$
Therefore, $\pi(P_v)=\sum\limits_{e:s(e)=v}\pi(S_e)\pi(S_e)^*$.
\fim  

The above theorem says that for a given E-branching system there exists a representation of $C^*(E)$ in $B(L^2(X,\mu))$. But that would be meaningless if E-branching systems did not exist. In the next section we show that this is not the case.

\section{Representations in $L^2(\R)$}

Next we show that for any given graph $E$, with $E^0$ and $E^1$ countable, there exists an E-branching system in $\R$ associated. We then show that for graphs that satisfy condition $(L)$ the representations arising from these E-branching systems in $\R$ are faithful. Our proof is constructive and one can actually obtain a great number of E-branching systems following the ideas below.

\begin{teorema}\label{existencebrancsys}
Let $E=(E^0,E^1,r,s)$ be a graph, with $E^0,E^1$ both countable. Then there exists an E-branching  system $(X,\mu)$, where the space $X$ is an (possible unlimited) interval of $\R$ and $\mu$ is the Lebesgue measure.
\end{teorema}

\demo 
Let $E^1=\{e_i\}_{i=1}^\infty$ (or, if $E^1$ is finite, let $E^1=\{e_i\}_{i=1}^N$). For each $i\geq 1$ define $R_{e_i}=[i-1,i]$. 
Let $W=\{v\in E^0\,:\,\,\text{ v is a sink}\}$ (a vertex $v\in E^0$ is a sink if $v\notin s(E^1)$). Note that $W$ is finite or infinite countable. Write $W=\{v_i\,:\,\,i=1,2,3,...\}$. For each $v_i\in W$, define $D_{v_i}=[-i,-i+1]$. For the vertices $u\in E^0$ which are not sinks, define $D_u=\bigcup\limits_{e_i:s(e_i)=u}R_{e_i}$. Note that items 1-4 from  definition \ref{brancsystem} are satisfied, considering the Lebesgue measure $\mu$. It remains to define functions which satisfy items 5-6.

Let $\overline{e}\in E^1$. 

If $r(\overline{e})$ is a sink then $r(\overline{e})=v_i\in W$, and so $D_{r(\overline{e})}=[-i, -i+1]$. Then we define $f_{\overline{e}}:D_{r(\overline{e})}\rightarrow R_{\overline{e}}$ as being a $C^1$-diffeomorphism (for example, the linear diffeomorphism). Note that such $f_{\overline{e}}$ in fact exists, because $D_{r(\overline{e})}$ and $R_{\overline{e}}$ are both closed limited intervals of $\R$. 

 If $r(\overline{e})=\overline{v}$ is not a sink, then $$D_{r(\overline{e})}=D_{\overline{v}}=\bigcup\limits_{e:s(e)=\overline{v}}R_e.$$ 
To define the function $f_{\overline{e}}:D_{r(\overline{e})}\rightarrow R_{\overline{e}}$ in this case we proceed as follows. 
 
 First we divide the interval $\stackrel{\circ}{R_{\overline{e}}}$ (where $\stackrel{\circ}{R_{\overline{e}}}$ denotes the interior of $R_{\overline{e}}$) in $\#\{e:s(e)=\overline{v}\}$ intervals $I_e$ (notice that we might have to divide $\stackrel{\circ}{R_{\overline{e}}}$ in a countable infinite number of intervals). Then, we define $\tilde{f_{\overline{e}}}:\bigcup\limits_{e:s(e)=\overline{v}}\stackrel{\circ}{R_e}\rightarrow \bigcup\limits_{e:s(e)=\overline{v}}\stackrel{\circ}{I_e}$ so that $\tilde{f_{\overline{e}}}_{|_{\stackrel{\circ}{R_e}}}$ is a $C^1$-diffeomorphism between $\stackrel{\circ}{R_e}$ and $\stackrel{\circ}{I_e}$ (for example, the linear diffeomorphism). We now define $f_{\overline{e}}:D_{r(\overline{e})}\rightarrow R_e$ as being a extension of $\tilde{f_{\overline{e}}}$ to $D_{r(\overline{e})}$ and $f_{\overline{e}}^{-1}:R_{\overline{e}}\rightarrow D_{r(\overline{e})}$ as being a extension of $\tilde{f_{\overline{e}}}^{-1}$ to $R_{\overline{e}}$.

For a given $e\in E^1$, the maps $f_e$ and $f_e^{-1}$ are measurable maps, in the measure space $D_{r(e)}$ and $R_e$, respectively, with the Lebesgue measure $\mu$ (and the Borel sets). Moreover, $\mu\circ f_e$ and $\mu \circ f_e^{-1}$ are $\sigma$-finite measures of (the measurable subsets of) $D_{r(e)}$ and $R_e$, respectively. It remains to see that there exists the nonnegative Radon-Nikodym derivatives $\Phi_{f_e}$ and $\Phi_{f_e}^{-1}$, and this follows from \cite{hewitt}.

Now, defining $$X=\left(\bigcup\limits_{e_i\in E^1}R_{e_i}\right)\cup\left(\bigcup\limits_{v_i\in W}D_{v_i}\right)$$
 we obtain the desired $(X,\mu)$ E-branching system.
 
\fim

Theorem \ref{existencebrancsys} guarantees that every graph C*-algebra (from a countable graph) may be represented in $B(L^2(\R))$. Of course when the graph C*-algebra is simple, the representations obtained via theorems \ref{rep} and \ref{existencebrancsys} are faithful.

Another case when the representations obtained via theorem \ref{rep} are faithful, without $C^*(E)$ being simple, is the case when the graph satisfies a special condition, called condition $(L)$. A graph $E$ satisfies condition $(L)$ if each loop has an exit, that is, if $x_1...x_n$ is a loop then there is a vertex $e$ such that $s(e)=s(x_i)$ for some $i$ but $e\neq x_i$. 

\begin{teorema} \label{conditionL} Let $E$ be a countable graph which satisfies the condition $(L)$. Then the representation $\pi:C^*(E)\rightarrow B(L^2(X))$, where $X$ is an (possible unlimited) interval of $\R$,  obtained via theorems \ref{existencebrancsys} and \ref{rep} is faithful.
\end{teorema}

\demo First note that for each $v\in E^0$, $P_v\in C^*(E)$ is a non-null element, because, by \ref{rep} and \ref{existencebrancsys}, there exists a representation $\pi:C^*(E)\rightarrow B(L^2(\mathbb{R}))$ such that $\pi(P_v)$ is the multiplication operator by $\chi_{D_v}$, where $D_v$ is a set of positive Lebesgue measure.  

Let $\pi:C^*(E)\rightarrow B(L^2(\R))$ be the representation obtained via \ref{existencebrancsys} and \ref{rep}. Since $E$ satisfies condition $(L)$,  by  [\cite{flr}:2], $\pi$ is faithful. 
\fim

\section{Nonsingular E-branching systems and Perron-Frobenius operators}\label{perronfrobenius}

Nonsingular maps\footnote{by a nonsingular map $F:X\rightarrow X$ we mean a measurable map such that $\mu(F^{-1}(A))=0$ if $\mu(A)=0$.} on a measure space, $(X,\mu)$, are of great interest in ergodic theory. In particular, each nonsingular map gives rise to a Perron-Frobenius operator on $L^1(X,\mu)$. In this section, we give a nice description of the Perron-Frobenius operator (for a large class of functions in $L^1(X,\mu)$) in terms of the representations introduced in the previous section. Unfortunately we can not do this for all nonsingular maps, but we can do it for all nonsingular maps that arise naturally from an E-branching system. 

To see how a nonsingular map arise from an E-branching system, recall that in the proof of theorem \ref{existencebrancsys} the measure space $(X,\mu)$ could be written as a disjoint union of the subsets $\{R_e\}_{e\in E^1}$ and $\{D_u\}_{u\in W}$ almost everywhere.
In this case, we may define a map $F:X\rightarrow X$ as follows:
$$F(x)=\left\{\begin{array}{cl} f_e^{-1}(x) & \text{ if } x\in \stackrel{o}{R_e}\text{ for some e }\in E^1 \\
0 & \text{ if }x\text{ is a extreme point of some interval } R_e \\
x & \text{ if }x \in D_u \text{ for some }u\in W
\end{array}\right .$$

The map $F$ above is nonsingular, since if $A\subseteq X$ is a measurable subset with $\mu(A)=0$ then, for a given $u\in W$, $\mu(F^{-1}(A)\cap D_u)=\mu(A\cap D_u)=0$. Furthermore, for each $e\in E^1$, $\mu(F^{-1}(A)\cap R_{e})=\mu(f_e(A\cap D_{r(e)}))=0$, since $\mu(A\cap D_{r(e)})=0$ and $\mu\circ f_e$ is absolutely continuous with respect to $\mu$ in $D_{r(e)}$.
So $\mu(F^{-1}(A))=0$, because $X$ is a countable union of the sets $R_e$ and $D_u$. Note also that the map $F$ defined above has the property that $F_{|_{R_e}}\stackrel{\mu-a.e.}{=}f_e^{-1}$. This motivates the definition of nonsingular E-branching systems:

\begin{definicao}
A nonsingular E-branching system $(X,\mu, F)$, associated to a directed graph $E=(E^0,E^1,r,s)$, is an E-branching system $(X,\mu)$, as defined in $\ref{brancsystem}$, together with a nonsingular map $F:X\rightarrow X$ such that $F_{|_{R_e}}\stackrel{\mu-a.e}{=}f_e^{-1}$ for each $e\in E^1$.
\end{definicao}

From theorem \ref{existencebrancsys} and the discussion above, we obtain promptly the corollary below.

\begin{corolario} For a given directed graph $E=(E^0,E^1,r,s)$, with $E^0,E^1$ both countable, there exists a nonsingular E-branching system $(X,\mu)$, where $X$ is a (possible unlimited) closed interval of $\R$ and $\mu$ the Lebesgue measure.
\end{corolario}

Actually, every E-branching system of a countable directed graph is also a nonsingular E-branching system, as we show below.

\begin{proposicao}\label{nonsingbranchsys} For every E-branching system $(X,\mu)$ of a countable directed graph $E=(E^0,E^1,r,s)$, there exists a map $F:X\rightarrow X$ such that $(X,\mu,F)$ is a nonsingular E-branching system.
\end{proposicao}

\demo Let $Y$ be the complement of $\bigcup\limits_{e\in E^1}R_e$ in $X$. Then we can write $X=Y\bigcup\limits^.(\bigcup\limits_{e\in E^1}R_e)$. Notice that, for each $e\in E^1$, there exists $\tilde{R_e}\subseteq R_e$ such that $\tilde{R_e}\stackrel{\mu-a.e.}{=}R_e$ and since $E^1$ is countable $\tilde{R_e}$ can be chosen such that all $\tilde{R_e}$ are pairwise disjoint. So, $X\stackrel{\mu-a.e.}{=}Y\bigcup\limits^.(\bigcup\limits_{e\in E^1}^.\tilde{R_e})$. Now, let $$\tilde{F}:Y\bigcup\limits^.(\bigcup\limits_{e\in E^1}^.\tilde{R_e})\rightarrow Y\bigcup\limits^.(\bigcup\limits_{e\in E^1}^.\tilde{R_e})$$
be defined by 
$$\tilde{F}(x)=\left\{\begin{array}{cl}f_e^{-1}(x) & \text{ if } x\in \tilde{R_e} \text{ for some }e\in E^1\\
x & \text{ if } x\in Y,
\end{array}\right.$$ and define $F:X\rightarrow X$ as an extension of $\tilde{F}$.

So, $F_{|_{R_e}}\stackrel{\mu-a.e}{=}f_e^{-1}$ for each $e\in E^1$. It remains to check that $F$ is nonsingular. Let $A\subseteq X$ be a measurable set with $\mu(A)=0$. Then, since $E^1$ is countable, it is enough to show that $\mu(F^{-1}(A)\cap Y)=0$ and $\mu(F^{-1}(A)\cap \tilde{R_e})=0$ for each $e\in E^1$. Now, note that $\mu(F^{-1}(A)\cap Y)=\mu(A\cap Y)=0$ and for each $e\in E^1$, $\mu(F^{-1}(A)\cap \tilde{R_e})=\mu(f_e(A\cap D_{r(e)}))=\mu\circ f_e(A\cap D_{r(e)})=0$, since $\mu\circ f_e$ is absolutely continuous with respect to $\mu$, in $D_{r(e)}$.
\fim

Recall that if $F:X\rightarrow X$ is a nonsingular map, where $X$ is a measure space with measure $\mu$, the Perron-Frobenius operator induced by $F$, denoted by $P_F$, is the operator in $B(L^1(X,\mu))$ such that for all $\psi\in L^1(X,\mu)$, and for all measurable subset $A\subseteq X$, the equality
$$\int\limits_AP_F(\psi)(x)d\mu=\int\limits_{F^{-1}(A)}\psi(x)d\mu$$ holds.

If $(X,\mu,F)$ is a nonsingular E-branching system of a directed graph $E$ then, by theorem \ref{rep}, there exists a *-representation of the graph algebra $C^*(E)$ in $B(L^2(X,\mu))$. The next theorem shows a relation between this representation and the Perron-Frobenius operator $P_F$.

\begin{teorema}\label{perronrep}
Let $E$ be a countable directed graph, $(X,\mu,F)$ be a nonsingular $E$-branching system and let $\varphi\in L^2(X,\mu)$.
\begin{enumerate}
\item If $supp(\varphi)\subseteq \bigcup\limits_{i=1}^N R_{e_i}$ then $P_F(\varphi^2)=\sum\limits_{i=1}^N\left(\pi(S_{e_i}^*)\varphi\right)^2.$

\item If $\varphi$ is a real valued function and $supp(\varphi)\subseteq \bigcup\limits_{i=1}^\infty R_{e_i}$ then $P_F(\varphi^2)=\lim\limits_{N\rightarrow \infty}\sum\limits_{i=1}^N\left(\pi(S_{e_i}^*)\varphi\right)^2,$ where the convergence occurs in the norm of $L^1(X,\mu)$.

\item If $supp(\varphi)\subseteq \bigcup\limits_{i=1}^\infty R_{e_j}$ and $\varphi=u+iv$, with $u,v$ real functions such that $uv\in L^2(X,\mu)$, then $P_F(\varphi^2)=\lim\limits_{N\rightarrow \infty}\sum\limits_{i=1}^N\left(\pi(S_{e_i}^*)\varphi\right)^2,$ where the convergence occurs in the norm of $L^1(X,\mu)$.
\end{enumerate}
\end{teorema}

\demo The first assertion follows from the fact that for each measurable set $A\subseteq X$, $\int\limits_AP_F(\varphi^2)(x)d\mu=\int\limits_A\sum\limits_{i=1}^N\left(\pi(S_{e_i}^*)\varphi(x)\right)^2d\mu$. We prove this equality below, and to do it we use the Radon-Nikodym derivative of $\mu\circ f_i$, the change of variable theorem and the fact that $F^{-1}(A)\cap R_{e_i}=f_{e_i}(A\cap D_i)$. 

Given a measurable set $A\subseteq X$, notice that
$$\sum\limits_{i=1}^N\int\limits_A\left(\pi(S_{e_i}^*)\varphi(x)\right)^2d\mu=\sum\limits_{i=1}^N\int\limits_A\chi_{D_{e_i}}(x)\Phi_{f_{e_i}}(x)\varphi(f_{e_i}(x))^2d\mu=$$
$$=\sum\limits_{i=1}^N\int\limits_{A\cap D_{e_i}}\Phi_{f_{e_i}}(x)\varphi(f_{e_i}(x))^2d\mu=\sum\limits_{i=1}^N\int\limits_{A\cap D_{e_i}}\varphi(f_{e_i}(x))^2d(\mu\circ f_{e_i})=$$
$$=\sum\limits_{i=1}^N\int\limits_{f_i(A\cap D_{e_i})}\varphi(x)^2d\mu=\sum\limits_{i=1}^N\int\limits_{F^{-1}(A)\cap R_{e_i}}\varphi(x)^2d\mu=$$ $$=\sum\limits_{i=1}^N\int\limits_{F^{-1}(A)}\chi_{R_{e_i}}\varphi(x)^2d\mu=\int\limits_{F^{-1}(A)}\sum\limits_{i=1}^N\chi_{R_{e_i}}\varphi(x)^2d\mu=$$
$$=\int\limits_{F^{-1}(A)}\varphi(x)^2d\mu=\int\limits_AP_F(\varphi^2)(x)d\mu.$$

Before we proceed with the proof of the second and third statements of the theorem, let us prove the following claim.\newline
Claim: If $h\subseteq L^1(X,\mu)$ is a real function with $supp(h)\subseteq\bigcup\limits_{j=1}^ \infty R_{e_j}$ then $\lim\limits_{N\rightarrow \infty}P_F(h_N)=P_F(h)$, where $h_N=\sum\limits_{j=1}^N\chi_{R_{e_j}}h$.

Suppose first that $h(x)\geq 0$ $\mu-a.e.$. Then, $(h_N)_{n\in \N}$ is a increasing sequence, bounded above by $h$, and so, 

$$\lim\limits_{N\rightarrow \infty}\int_XP_F(h_N)(x)d\mu=\lim\limits_{N\rightarrow \infty}\int_X h_N(x)d\mu=$$ $$=\int_X h(x) d\mu=\int_X P_F(h)(x) d\mu.$$ Now, since $h_N\leq h$ we have that $P_F(h_N)\leq P_F(h)$ and hence 
$$\lim\limits_{N\rightarrow \infty}\|P_F(h)-P_F(h_N)\|_1=\lim\limits_{N\rightarrow \infty}\int\limits_X |P_F(h)(x)-P_F(h_N)(x)|d\mu=$$
$$=\lim\limits_{N\rightarrow \infty}\int\limits_X P_F(h)(x)-P_F(h_N)(x)d\mu=0.$$

To prove the claim for a real function $h\in L^1(X,\mu)$, write $h=h_1-h_2$, where $h_1$ and $h_2$ are nonnegative functions and use the linearity of $P_F$.

Next we prove the second statement of the theorem. Define $\varphi_N=\sum\limits_{j=1}^NR_{e_j}\varphi$. By the above claim, $\lim\limits_{N\rightarrow \infty}P_F(\varphi_N^2)=P_F(\varphi^2)$. By the first statement, $\sum\limits_{j=1}^N(\pi(S_{e_j})^*\varphi_N)^2=P_F(\varphi_N^2)$, and a simple calculation shows that
$$\sum\limits_{i=1}^N\left(\pi(S_{e_i}^*)\varphi_N\right)^2=\sum\limits_{i=1}^N\left(\pi(S_{e_i}^*)\varphi\right)^2.$$

So, we conclude that $$\lim\limits_{N\rightarrow \infty}\sum\limits_{i=1}^N\left(\pi(S_{e_i}^*)\varphi\right)^2=P_F(\varphi^2).$$

To prove the third statement, let $\varphi\in L^2(X,\mu)$ be a complex function and write $\varphi=u+iv$, with $u,v$ real functions. Define $u_N=\sum\limits_{j=1}^N\chi_{R_{e_j}}.u$, and $v_N=\sum\limits_{j=1}^N\chi_{R_{e_j}}.v$. Then, 
$$\sum\limits_{j=1}^N(\pi(S_{e_j})^*\varphi)^2=\sum\limits_{j=1}^N(\pi(S_{e_j})^*(u_N+iv_N))^2=$$ 
$$=\sum\limits_{j=1}^N(\pi(S_{e_j})^*u_N)^2-\sum\limits_{j=1}^N(\pi(S_{e_j})^*v_N)^2+i\sum\limits_{j=1}^N 2\chi_{D_{e_j}}.\Phi_{f_{e_j}}.(u_N\circ f_{e_j}).(v_N\circ f_{e_j})=$$ 

$$=P_F(u_N^2)-P_F(v_N^2)-i2P_F(u_Nv_N).$$

The last equality follows from the first statement of the theorem and from the fact that 
$$\chi_{D_{e_j}}.\Phi_{f_{e_j}}.(u_N\circ f_{e_j}).(v_N\circ f_{e_j})=P_F(u_Nv_N),$$ since
for each $E\subseteq X$, 
$$\int\limits_E\chi_{D_{e_j}}(x)\Phi_{f_{e_j}}(x)(u_N(f_{e_j}(x)))(v_N(f_{e_j}(x)))d\mu=\int\limits_EP_F(u_Nv_N)d\mu.$$  

Finally, since $u_N,v_N$ are real functions, by the Claim proved above,
$$\lim\limits_{N\rightarrow \infty}\sum\limits_{j=1}^N(\pi(S_{e_j})^*\varphi)^2=\lim\limits_{N\rightarrow \infty}P_F(u_N^2)-P_F(v_N^2)-i2P_F(u_Nv_N)=$$
$$=P_F(u^2)-P_F(v^2)+i2P_F(uv)=P_F((u+iv)^2)=P_F(\varphi^2).$$

\fim

\section{Examples}

We finish the paper with two examples of how our construction works.

\begin{exemplo}\label{compact} The compact operators in a separable Hilbert space. 
\end{exemplo}

Consider the following graph $E$

\centerline{
\setlength{\unitlength}{2cm}
\begin{picture}(4,0.6)
\put(0,0){$\dots$}
\put(0.5,0){\circle*{0.08}}
\put(0.6,0){\line(1,0){1}}
\put(1.1,-0.05){$>$}
\put(0.4,0.1){$v_{-1}$}
\put(1.1,0.1){$e_0$}
\put(1.7,0){\circle*{0.08}}
\put(1.8,0){\line(1,0){1}}
\put(1.6,0.1){$v_0$}
\put(2.2,0.1){$e_1$}
\put(2.25,-0.05){$>$}
\put(2.9,0){\circle*{0.08}}
\put(2.9,0){\line(1,0){1}}
\put(2.8,0.1){$v_1$}
\put(3.4,0.1){$e_2$}
\put(3.4,-0.05){$>$}
\put(4,0){\circle*{0.08}}
\put(3.9,0.1){$v_2$}
\put(4.1,0){\dots}
\end{picture}}
\vspace{0.6cm}
The graph C*-algebra, $C^*(E)$, is the algebra of compact operators in a separable Hilbert space, which we denote by $\mathcal{K}$. 

First we will show how to use our methods to obtain a faithful representation of $\mathcal{K}$ in $B(L^2(\R))$. We could follow the steps of the proof of theorem $\ref{existencebrancsys}$, but due to the simetry of this graph we will build an E-branching system in the following way:

Let $R_{e_i}=[i-1,i]$ and so we must define $D_{v_i}=[i,i+1]$ for each $i\in \Z$. Also, let $f_{e_i}:D_{r(e_i)}\rightarrow R_{e_i}$ be defined by $f_{e_i}(x)=x-1$. One can check that this defines an $E$-branching system $(\R,\mu)$, where $\mu$ is the Lebesgue measure in $\R$. Following theorem \ref{rep} we obtain a representation $\pi:C^*(E)\rightarrow B(L^2(\R))$, such that $\pi(P_{v_i})$ is the multiplication operator by $\chi_{[i,i+1]}$ and $\pi(S_{e_i})\phi_{|_x}=\chi_{[i-1,i]}(x)\phi(x+1)$ for each $\phi \in L^2(\R)$ and $x\in \R$. By theorem \ref{conditionL} this is a faithful representation. 

Following section \ref{perronfrobenius} (see proposition \ref{nonsingbranchsys} and theorem \ref{perronrep}), there is a Perron-Frobenius operator, $P_F$, associated to the $E$-branching system $(\R,\mu)$. The nonsingular map $F:\R\rightarrow \R$ restricted to each $R_{e_i}$ is the inverse of $f_{e_i}$ and  is given below 

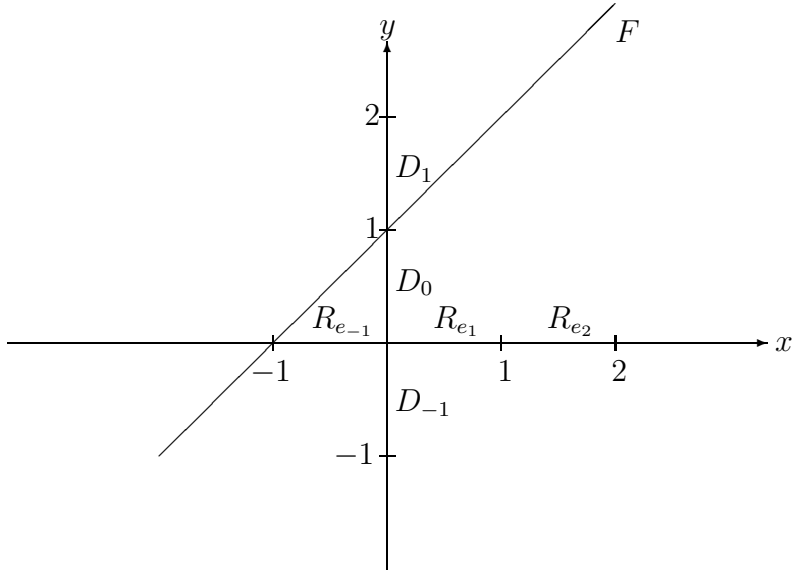
\setlength{\unitlength}{1cm}
\begin{picture}(5,5)
\put(0,0){\vector(1,0){10}}
\put(10.1,-0.1){$x$}
\put(5,-3){\vector(0,1){7}}
\put(4.9,4.1){$y$}
\put(6.5,-0.1){\line(0,1){0.2}}
\put(6.45,-0.5){$1$}
\put(5.6,0.2){$R_{e_1}$}
\put(8,-0.1){\line(0,1){0.2}}
\put(7.95,-0.5){$2$}
\put(7.1,0.2){$R_{e_2}$}
\put(3.5,-0.1){\line(0,1){0.2}}
\put(3.2,-0.5){$-1$}
\put(4,0.2){$R_{e_{-1}}$}
\put(4.9,-1.5){\line(1,0){0.2}}
\put(4.3,-1.6){$-1$}
\put(5.1,-0.9){$D_{-1}$}
\put(4.9,1.5){\line(1,0){0.2}}
\put(4.7,1.4){$1$}
\put(5.1,0.7){$D_0$}
\put(4.9,3){\line(1,0){0.2}}
\put(4.7,2.9){$2$}
\put(5.1,2.2){$D_1$}
\put(2,-1.5){\line(1,1){6}}
\put(8,4){$F$}

\end{picture}
\vspace{3cm}

The Perron-Frobenius operator in this case is easy to calculate, and is given by $P_F(\psi)=\psi\circ F^{-1}$ for each $\psi\in L^1(\R)$. In particular, for $\varphi\in L^2(\R)$ with $supp(\varphi)\subseteq \bigcup\limits_{i=1}^\infty R_{e_i}$, this characterization may also be obtained via theorem \ref{perronrep}.   

\begin{exemplo}
\end{exemplo}
Let $E$ be the finite graph below:

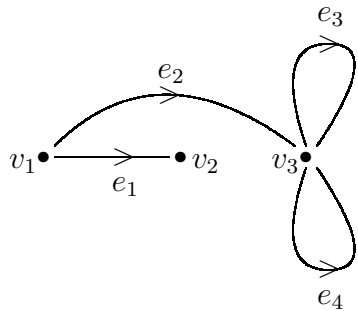
\centerline{
\setlength{\unitlength}{1.5cm}
\begin{picture}(3,1.5)
\put(0,0){\circle*{0.1}}
\put(0.1,0){\line(1,0){1}}
\put(0.6,-0.07){$>$}
\put(1.2,0){\circle*{0.1}}
\put(2.3,0){\circle*{0.1}}
\qbezier(0.1,0.1)(1,1)(2.2,0.1)
\put(1,0.475){$>$}
\qbezier(2.3,0.1)(2,1)(2.5,1)
\qbezier(2.4,0.1)(3,1)(2.5,1)
\put(2.4,0.925){$>$}
\qbezier(2.3,-0.1)(2,-1)(2.5,-1)
\qbezier(2.4,-0.1)(3,-1)(2.5,-1)
\put(2.4,-1.06){$>$}
\put(-0.3,-0.1){$v_1$}
\put(1.3,-0.1){$v_2$}
\put(2,-0.1){$v_3$}
\put(0.6,-0.3){$e_1$}
\put(1,0.7){$e_2$}
\put(2.4,1.2){$e_3$}
\put(2.4,-1.3){$e_4$}
\end{picture}}
\vspace{2cm}

Then $C^*(E)$ is a non-simple C*-algebra (see \cite{BHRS} for its ideal structure), but $E$ is a graph that satisfies condition $(L)$ and hence we can aply theorem \ref{conditionL} to obtain a faithful representation of $C^*(E)$ in $B(L^2(X))$, where $X$ arises from theorem \ref{existencebrancsys}. Therefore, to construct an E-branching system, we will follow the steps of the proof of theorem \ref{existencebrancsys}. 

So, let $R_{e_i}=[i-1,i]$ for $i\in \{1,2,3,4\}$, $D_{v_1}=[0,2]$, $D_{v_2}=[-1,0]$, $D_{v_3}=[2,4]$, and let $f_{e_i}:D_{r(e_i)}\rightarrow R_{e_i}$ be the affine maps defined as in the figure below:

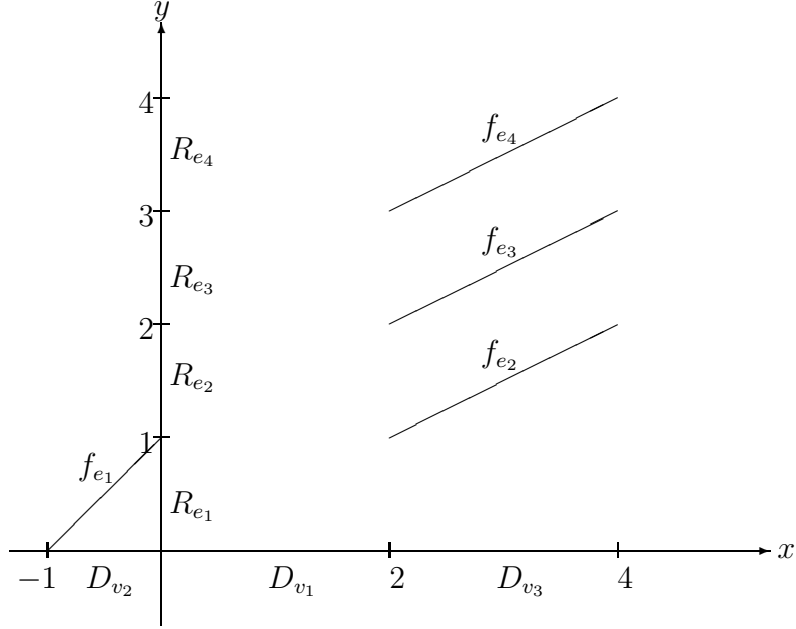
\centerline{
\setlength{\unitlength}{1cm}
\begin{picture}(10,8)
\put(0,0){\vector(1,0){10}}
\put(10.1,-0.1){$x$}
\put(2,-1){\vector(0,1){8}}
\put(1.9,7.1){$y$}
\put(0.5,-0.1){\line(0,1){0.2}}
\put(0.1,-0.5){$-1$}
\put(1,-0.5){$D_{v_2}$}
\put(5,-0.1){\line(0,1){0.2}}
\put(5,-0.5){$2$}
\put(3.4,-0.5){$D_{v_1}$}
\put(8,-0.1){\line(0,1){0.2}}
\put(8,-0.5){$4$}
\put(6.4,-0.5){$D_{v_3}$}
\put(1.9,1.5){\line(1,0){0.2}}
\put(1.7,1.3){$1$}
\put(2.1,0.5){$R_{e_1}$}
\put(1.9,3){\line(1,0){0.2}}
\put(1.7,2.8){$2$}
\put(2.1,2.2){$R_{e_2}$}
\put(1.9,4.5){\line(1,0){0.2}}
\put(1.7,4.3){$3$}
\put(2.1,3.5){$R_{e_3}$}
\put(1.9,6){\line(1,0){0.2}}
\put(1.7,5.8){$4$}
\put(2.1,5.2){$R_{e_4}$}
\put(0.5,0){\line(1,1){1.5}}
\put(0.9,1){$f_{e_1}$}
\put(5,1.5){\line(2,1){3}}
\put(6.2,2.5){$f_{e_2}$}
\put(5,3){\line(2,1){3}}
\put(6.2,4){$f_{e_3}$}
\put(5,4.5){\line(2,1){3}}
\put(6.2,5.5){$f_{e_4}$}
\end{picture}}

\vspace{1cm}
where, for example, the map $f_{e_2}:D_{r(e_2)}=[2,4]\rightarrow [1,2]=R_{e_2}$ is defined by $f_{e_2}(x)=\frac{x}{2}$ for each $x\in [2,4]$.

In this example, the measure space is the interval $[-1,4]$, with the Lebesgue measure. The representation $\pi:C^*(E)\rightarrow B(L^2([-1,4]))$ induced by this $E$-branching system is such that: $\pi(P_{v_i})$ is the multiplication operator by $\chi_{D_{v_i}}$ for $i\in\{1,2,3\}$, $\pi(S_{e_1})\varphi=\chi_{[0,1]}\cdot\varphi\circ f_{e_1}^{-1}$ and $\pi(S_{e_i})\varphi=\sqrt{2}\chi_{[i-1,i]}\cdot \varphi\circ f_{e_i}^{-1}$ for each $i\in \{2,3,4\}$ and for each $\varphi\in L^2([-1,4])$.

Following section \ref{perronfrobenius}, there is a nonsingular map $F:[-1,4]\rightarrow [-1,4]$ associated to this $E$-branching system. The graph of $F$ is shown in the following figure:

\centerline{
\setlength{\unitlength}{1cm}
\begin{picture}(10,8)
\put(0,0){\vector(1,0){10}}
\put(10.1,-0.1){$x$}
\put(2,-2){\vector(0,1){9}}
\put(1.9,7.1){$y$}
\put(1.9,-1.5){\line(1,0){0.2}}
\put(1.4,-1.6){$-1$}
\put(0.5,-0.1){\line(0,1){0.2}}
\put(0.1,-0.5){$-1$}
\put(5,-0.1){\line(0,1){0.2}}
\put(4.95,-0.5){$2$}
\put(3.5,-0.1){\line(0,1){0.2}}
\put(3.45,-0.5){$1$}
\put(6.5,-0.1){\line(0,1){0.2}}
\put(6.45,-0.5){$3$}
\put(8,-0.1){\line(0,1){0.2}}
\put(8,-0.5){$4$}
\put(1.9,1.5){\line(1,0){0.2}}
\put(1.7,1.3){$1$}
\put(1.9,3){\line(1,0){0.2}}
\put(1.7,2.8){$2$}
\put(1.9,4.5){\line(1,0){0.2}}
\put(1.7,4.3){$3$}
\put(1.9,6){\line(1,0){0.2}}
\put(1.7,5.8){$4$}
\put(0.5,-1.5){\line(1,1){1.5}}
\put(2.05,-1.45){\line(1,1){1.45}}
\put(2,-1.5){\circle{0.1}}
\put(3.5,3){\line(1,2){1.45}}
\put(3.47,2.95){\circle{0.1}}
\put(4.97,5.95){\circle{0.1}}
\put(5,3){\line(1,2){1.45}}
\put(4.97,2.95){\circle{0.1}}
\put(6.47,5.95){\circle{0.1}}
\put(6.5,3){\line(1,2){1.45}}
\put(6.47,2.95){\circle{0.1}}
\put(7.97,5.95){\circle{0.1}}
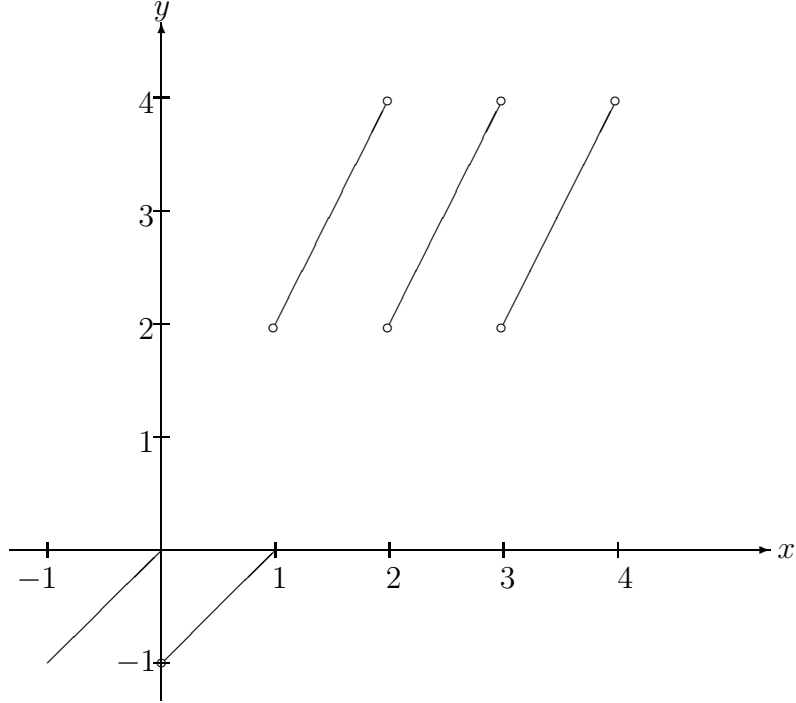
\put(1.5,-2.5){\text{the nonsingular map} $F:[-1,4]\rightarrow [-1,4]$}
\end{picture}}
\vspace{3.5cm}

This nonsingular map induces a Perron-Frobenius operator $P_F:L^1([-1,4])\rightarrow L^1([-1,4])$, and following theorem \ref{perronrep}, for each $\varphi\in L^2([-1,4])$ with $\supp(\varphi)\subseteq [0,4]$, we have that 

$$P_F(\varphi^2)=\sum\limits_{i=1}^4(\pi(S_e^*)\varphi)^2=\chi_{[-1,0]}\cdot(\varphi\circ f_{e_1})^2+\frac{1}{2}\chi_{[2,4]}\cdot\left[(\varphi\circ f_{e_2})^2+(\varphi\circ f_{e_3})^2+(\varphi\circ f_{e_4})^2\right].$$ So, for each $x\in [-1,4]$, it holds that

$$P_F(\varphi^2)_{|_x}=\chi_{[-1,0]}(x)\cdot\varphi(f_{e_1}(x))^2+\frac{1}{2}\chi_{[2,4]}(x)\cdot\left[\varphi(f_{e_2}(x))^2+\varphi(f_{e_3}(x))^2+\varphi(f_{e_4}(x))^2\right]=$$

$$=\chi_{[-1,0]}(x)\cdot\varphi(x+1)^2+\frac{1}{2}\chi_{[2,4]}(x)\cdot\left[\varphi\left(\frac{x}{2}\right)^2+\varphi\left(\frac{x}{2}+1\right)^2+\varphi\left(\frac{x}{2}+2\right)^2\right],$$
and this explicitly describes the Perron-Frobenius operator $P_F$ for a large number of functions.

\bibliographystyle{amsplain}
\bibliography{newbib}

\vspace{1.5pc}
\begin{center}
D. Gonçalves (daemig@gmail.com) and D. Royer (royer@mtm.ufsc.br)\\

Departamento de Matemática - Universidade Federal de Santa Catarina, Florianópolis, 88040-900, Brazil
\end{center}

\end{document}